# REPLENISHMENT COST OF A STORAGE SYSTEM MODELLED THROUGH A DEFECTIVE RENEWAL EQUATION

Manuel Alberto M. Ferreira

manuel.ferreira@iscte.pt

Instituto Universitário de Lisboa (ISCTE-IUL), BRU-IUL, ISTAR-IUL, Lisboa,

PORTUGAL

**ABSTRACT**

A storage system where a *k* units stock is replaced instantaneously, when required, is presented in this work. It is also supposed a Poisson demand. It is proved that the storage system replenishment cost expected present value function fulfills a defective renewal equation, being proposed for it an asymptotic expression.

## INTRODUCTION

For a certain item suppose that demand follows a Poisson distribution with parameter $\mu$:

$$P[\text{Demand up to time } t = x] = \frac{e^{-\mu t}(\mu t)^x}{x!}, x = 0,1,\ldots; t > 0.$$

So the time length availability of a *k* units stock, call it *X*, is distributed in accordance with a gamma probability distribution:

$$F(x) = \int_0^x f(y)dy\,;\ f(x) = \frac{\mu^k x^{k-1} e^{-\mu x}}{(k-1)!}, x > 0 \qquad (1.1)$$

with Laplace transform

$$\phi^k(s) = \left(\frac{\mu}{s+\mu}\right)^k, s > -\mu \qquad (1.2).$$

Suppose the following protocol for the storage system. At time $t = 0$ there are *k* units available in stock. Every time the system goes out of the stock, an amount of *k* units is instantaneously replaced.

The system will be evaluated through the stochastic process $W_k(t)$ and the random variable $V_k(t)$:

$$W_k(t) = \sum_{n=1}^{N(t)} \theta_k \prod_{m=1}^{n} e^{-rX_m}, W_k(t) = 0 \text{ if } N(t) = 0 \quad (1.3)$$
$$V_k = W_k(\infty) \quad (1.4)$$

where $W_k(t)$ is the stock replenishment cost present value, up to time $t$, when a payment of $\theta_k$ is made each time k units are replaced. The parameter $r, r > 0$ represents a deterministic discount interest rate. The collection $X_1, X_2, \ldots$ is a sequence of *iid* random variables with the same distribution as $X$. They represent the time intervals between successive stock replacements. The associated counting process $N(t)$ is defined as:

$$N(t) = \sup\{n : \sum_{i=1}^{n} X_i \leq t\}, N(t) = 0 \text{ if } X_1 > t$$

and $V_k$ is the storage system perpetual cost.

A particular attention will be given to (1.3) and (1.4) expected values:

$$w_k(t) = E[W_k(t)]; \ v_k = E[V_k] \quad (1.5).$$

**Notes:**

-Alternatively, instead of $\theta_k$, the cost of $k$ units, a time dependent value with exponential growth: $\theta_k = \theta_k(t) = \tilde{\theta}_k e^{\tilde{r}t}$ can be considered. In this case, the present value process of the stock replenishment cost would assume the form

$$W_k(t) = \sum_{n=1}^{N(t)} \tilde{\theta}_k \prod_{m=1}^{n} e^{-(r-\tilde{r})X_m}, W_k(t) = 0 \text{ if } N(t) = 0.$$

With this assumption the analysis performed for (1.3) and (1.4) still stands for $\tilde{r} < r$, being enough to replace $r$ by $(r - \tilde{r})$ in the above formulas.

-It is possible to consider $\theta_k$ as the acquisition and commercialization of $k$ items global result. For instance $\theta_k = bk - a$ where $a$ is the fixed cost associated to any operation concerning the stocks reposition and $b$ the unitary margin associated with buying and selling one item. Evidently, $b > \frac{a}{k}$.

**RESULTS**

**Proposition 2.1**

$$v_k = \frac{\theta_k \phi^k(r)}{1 - \phi^k(r)} \quad (2.1).$$

**Dem.:** Note that $V_k$ is a random perpetuity, see[8] and [10]. So it may be written as the solution of the random equation $V_k \stackrel{d}{=} e^{-rX}(\theta_k + V_k)$, $X$ and $V_k$ independent, being the identity in probability distribution. Applying expectations in this equation both sides it results (2.1). ∎

**Proposition 2.2**

The expected value function $w_k(t)$ satisfies the defective renewal equation[1]

$$w_k(t) = \theta_k \phi^k(r) F(t) + \int_0^t w_k(t-s) \phi^k(r) f(s) ds \quad (2.2).$$

**Dem.:** Conditioning to $N(t)$ and applying the expectation tower property,

$$w_k(t) = E[E[W_k(t)|N(t)]] = \theta_k \phi^k(r) \frac{1 - \gamma(t, \phi^k(r))}{1 - \phi^k(r)} \quad (2.3),$$

where $\gamma(t, s)$ is the $N(t)$ probability generating function. Call $F^{*n}$ the $n$-fold convolution of $F$ with itself and assume $F^{*0} = 1$ for all $t \geq 0$. Then

$$\gamma(t, s) = E[s^{N(t)}] = \sum_{n=0}^{\infty} s^n (F^{*n}(t) - F^{*(n+1)}(t)) = 1 - (1-s) \sum_{n=0}^{\infty} s^{n-1} F^{*n}(t).$$

Substituting in (2.3): $w_k(t) = \theta_k \phi^k(r) \sum_{n=1}^{\infty} \phi^{k(n-1)}(r) F^{*n}(t)$ and applying Laplace transforms in this equality both sides,

$$\int_0^{\infty} e^{-st} w_k(t) dt = \frac{\theta_k \phi^k(r) \phi^k(s)}{s(1 - \phi^k(r) \phi^k(s))} \quad (2.4)$$

which after rearranging and inversion returns (2.2). ∎

Defining $J(t) = w_k(\infty) - w_k(t)$, and noting that $w_k(\infty) = v_k$, let $j(t) = \frac{\theta_k \phi^k(r)}{1 - \phi^k(r)} (1 - F(t))$. Consider also the positive constant $\rho$ that satisfies $\int_0^{\infty} e^{\rho s} \phi^k(r) f(s) ds = \phi^k(r) \phi^k(-\rho) = 1$. That is $\rho = \frac{r\mu}{r+\mu}$, recalling (1.2).

**Proposition 2.3**

For the function $e^{\rho t} J(t)$ the following asymptotic expansion is true:

$$\lim_{t \to \infty} e^{\rho t} J(t) = \frac{1}{\mu_0} \int_0^{\infty} e^{\rho s} j(s) ds \quad (2.5),$$

$$\mu_0 = \int_0^{\infty} s e^{\rho s} \phi^k(r) f(s) ds \quad (2.6).$$

**Dem.:** Computing $J(t)$ from (2.2) and multiplying the obtained equation both members by $e^{\rho t}$, the following common renewal equation is achieved:

$$e^{\rho t} J(t) = e^{\rho t} j(t) + \int_0^t e^{\rho(t-s)} J(t-s) e^{\rho s} \phi^k(r) f(s) ds.$$

---

[1] See [9]

This proposition is a consequence of the key renewal theorem, see [9] and [1], application to this equation. ∎

**Notes:**

- The conclusions in the former propositions remain true - with eventually minor adaptations - for a lot of $X$ continuous probability distributions, if they have at least the nonnegative reals set as support.

-To specify those results for the special case of the initially defined gamma distribution, make repeated use of the Laplace transform (1.2) and note that the (2.5) right side can be written as: $\frac{1}{\mu_0}\int_0^\infty e^{\rho s}j(s)ds = \frac{1}{\frac{k(r+\mu)}{\mu^2}} \cdot \frac{\theta_k(r+\mu)}{r\mu} = \frac{\theta_k \mu}{kr}$. Using this result, Proposition 2.1 and Proposition 2.3 can be summarized as

$$v_k = \frac{\theta_k}{\alpha^k - 1}, \quad \alpha = \frac{r}{\mu} + 1 \quad (2.7),$$

$$w_k(t) \approx \frac{\theta_k}{\alpha^k - 1} - \frac{\theta_k}{k(\alpha-1)} e^{-\frac{r}{\alpha}t} \quad (2.8).$$

-The approximation in (2.8) must be faced with the same meaning as in (2.7). For more detail, in these results, see [2].

### EXAMPLES

### Example 3.1

Consider the case $k = 1$. Then (2.8) is $w_1(t) = \frac{\theta_1}{\alpha-1}\left(1 - e^{-\frac{r}{\alpha}t}\right)$. The right side Laplace transform is $\int_0^\infty e^{-st}w_1(t)dt = \frac{\theta_1 r}{s(s\alpha+r)(\alpha-1)}$. Calculating the (2.4) right side the same result is got. Indeed, in this case, the given expression is exact. ∎

### Example 3.2

Now consider, as suggested above, $\theta_k = bk - a$. For the situation: $\mu = a = b = 1$ and $r = 0.02$, (2.7) is $v_k = \frac{k-1}{1.02^k - 1}$. The most interesting value for $k$ is so $\max_k\{v_k\} = 10$. In this situation it is obtained: $w_{10}(t) \approx \frac{9}{1.02^{10}-1} - 45e^{-\frac{t}{101}}$.

**Table 3.1**

| t | 10 | 20 | 50 | 100 | 200 | 500 | ∞ |
|---|---|---|---|---|---|---|---|
| $w_{10}(t)$ | .339 | 4.181 | 13.688 | 24.378 | 34.885 | 40.778 | 41.097 |

In Table 3.1 some values for this function are presented. ∎